\documentclass[12pt]{article}
\usepackage{amsfonts, amsthm, amsmath}

\def\FF{\mathbb{F}}
\def\PP{\mathbb{P}}
\def\QQ{\mathbb{Q}}

\def\Spec{{\mathrm{Spec\/}}}
\def\cris{{\mathrm{crys}}}

\begin{document}

\title{Crew's Euler Characteristic Formula Fails for Nonzero Slopes}
\author{Kiran S. Kedlaya}
\date{May 15, 2001}

\maketitle


Let $f: X \to Y$ be a finite \'etale cover of separated schemes of
fnite type over a field $k$ of characteristic
$p>0$. If $l\neq p$,
the Euler-Poincar\'e formula (a special case of Riemann-Hurwitz)
states that $\chi_c(X, \QQ_l) =
(\deg f)\chi(Y, \QQ_l)$,
where $\chi_c$ denotes the Euler characteristic
\[
\chi_c(X, \QQ_l) = \sum_{i=0}^{2\dim X} \dim_{\QQ_l} H_c^i(X, \QQ_l)
\]
with respect to \'etale cohomology with compact supports.

This relation fails in general for $l=p$, but holds in an important
special case. Suppose $f$ is a Galois cover and $\deg f$ is a power of $p$.
Then the relation $\chi_c(X, \QQ_p) = (\deg f)\chi(Y, \QQ_l)$ holds;
for complete curves, this is due to Shafarevich
\cite{shaf}, for arbitrary curves it is equivalent to the Deuring-Shafarevich
formula, and in general it is due to 
 in dimension 1 and by Crew \cite{crew} in general.

For $X$ and $Y$ smooth and proper over $\Spec\, k$, this assertion can
be reinterpreted in terms of crystalline cohomology. Namely, 
in this case, the crystalline cohomology groups $H_{\cris}^i(X/W)$ are
finitely generated modules over the Witt ring $W$ of $k$ with a
semilinear endomorphism $F$ (the Frobenius). For any rational number 
$\lambda$, we let $h_\lambda^i(X)$ be the multiplicity of the slope
$\lambda$ in $H_{\cris}^i(X/W)$; then $h_0^i(X) = \dim H_c^i(X, \QQ_l)$.
If we put
\[
\chi_\lambda(X) = \sum_{i=0}^{2\dim X} (-1)^i h_\lambda^i(X),
\]
then Crew's theorem is that $\chi_0(X) = |G| \chi_0(Y)$.

Crew asks whether $\chi_\lambda(X) = |G| \chi_\lambda(Y)$
for other values of $\lambda$. We present an explicit counterexample
to this assertion with $p=2$ and $\dim X = 1$.

Let $k = \overline{\FF_2}$. Consider the curves
\begin{align*}
C: u^2-u &= \frac{1 + x^2 + x^8 + x^{14} + x^{18}}{x^{21}} \\
D: v^2-v &= \frac{1}{x+1} \\
Y: w^2-w &= \frac{1 + x^2 + x^8 + x^{14} + x^{18}}{x^{21}} + \frac{1}{x+1}
\end{align*}
and let $X$ be the fibre product of $C$ and $D$ over the maps to $\PP_1$
given by $x$. Then $X$ admits a map to $Y$ given by setting $w=u+v$;
it is easily verified that this map is \'etale. Furthermore, using
Riemann-Hurwitz, one
calculates $g(C) = 10$, $g(D) = 0$, $g(Y) = 11$, and $g(X) = 21$.

The characteristic polynomials of Frobenius on $H^1(C), H^1(X), H_1(Y)$ can be
obtained by counting points over finite extensions of $\FF_2$. We compute
\begin{align*}
P_C(t) &= 1 - 32t^{10} + 1024t^{20} \\
P_Y(t) &= 1+t+2t^2+4t^3+4t^4+4t^5+8t^6+8t^7+8t^8+16t^9+32t^{10}
+32t^{11} \\
&+64t^{12}+64t^{13}+64t^{14}+128t^{15}+256t^{16}+256t^{17}+512t^{18}+1024t^{19}\\
&
+1024t^{20}+1024t^{21}+2048t^{22}
\end{align*}
and $P_X(t) = P_C(t)P_Y(t)$; of course $P_D(t) = 1$. From these polynomials
we can read off the $h_\lambda^i$; namely,
\begin{gather}
h_{0}^1(Y) = h_1^1(Y) = 1, h_{3/7}^1(Y) = h_{4/7}^1(Y) = 7, h_{1/2}^1(Y) = 6 \\
h_{0}^1(X) = h_1^1(X) = 1, h_{3/7}^1(X) = h_{4/7}^1(X) = 7, h_{1/2}^1(X) = 26
\end{gather}
and so $\chi_\lambda(X) \neq 2\chi_\lambda(Y)$ for $\lambda \in \{
3/7, 1/2, 4/7\}$.

The underlying phenomenon seems to be that while $C$ is supersingular,
the generic curve of the form $u^2-u = A(x)/x^{21}$ has slopes
$3/7$ and $4/7$ with multiplicity 7. It appears that the supersingularity
of $C$ is unstable under fiber products.

\section*{Acknowledgments}
Thanks to Arthur Ogus for directing the author to Crew's thesis,
to Richard Crew for historical remarks, and to Hui June Zhu for
numerical data that led to this example.
The author was supported by an NSF Postdoctoral Fellowship.


\begin{thebibliography}{9}

\bibitem{crew}
R. Crew, Etale $p$-covers in characteristic $p$, \textit{Comp. Math.}
\textbf{52} (1984), 31--45.

\bibitem{shaf}
I. Shafarevich, On $p$-extensions, \textit{Mat. Sbornik} \textbf{20}
(1947), 351--363.

\end{thebibliography}
\end{document}